\documentclass[conference]{IEEEtran}

\usepackage{amsmath,amssymb,amsthm}
\usepackage{graphicx}
\usepackage{bm,color}
\usepackage{algorithmicx,algorithm,algpseudocode}
\usepackage{url}
\usepackage{cite}
\usepackage{textcomp}

\allowdisplaybreaks[3]

\newtheorem{definition}{Definition}
\newtheorem{remark}{Remark}

\newcommand{\bp}{\mathbf{p}}
\newcommand{\bv}{\mathbf{v}}

\newcommand{\be}{\mathbf{e}}
\newcommand{\bl}{\bm{\ell}}

\newcommand{\bw}{\mathbf{w}}

\newcommand{\bu}{\mathbf{u}}
\newcommand{\bb}{\mathbf{b}}
\newcommand{\bz}{\mathbf{z}}

\newcommand{\pg}[2]{P_{G_{#1}}^{#2}}
\newcommand{\pd}[2]{P_{D_{#1}}^{#2}}
\newcommand{\pb}[2]{P_{B_{#1}}^{#2}}
\newcommand{\bj}[2]{B_{#1}^{#2}}
\newcommand{\bjmax}[1]{B^{\text{max}}_{#1}}
\newcommand{\bjmin}[1]{B^{\text{min}}_{#1}}
\newcommand{\pgmax}[1]{P^{\text{max}}_{{G}_{#1}}}
\newcommand{\pgmin}[1]{P^{\text{min}}_{{G}_{#1}}}
\newcommand{\pdmax}[1]{P^{\text{max}}_{{D}_{#1}}}
\newcommand{\pdmin}[1]{P^{\text{max}}_{{D}_{#1}}}
\newcommand{\pbmax}[1]{P^{\text{max}}_{{B}_{#1}}}
\newcommand{\pbmin}[1]{P^{\text{max}}_{{B}_{#1}}}
\newcommand{\rup}[1]{R^{\text{up}}_{#1}}
\newcommand{\rdn}[1]{R^{\text{dn}}_{#1}}
\newcommand{\Prob}{\mathbb{P}}
\newcommand{\R}{\mathbb{R}}
\newcommand{\1}{1{\hskip -2.55 pt}\hbox{I}}
\newcommand{\cK}{{\cal K}}
\newcommand{\cP}{{\cal P}}

\newcommand{\cM}{{\cal M}}
\newcommand{\cN}{{\cal N}}
\newcommand{\cI}{{\cal I}}

\newcommand{\cT}{{\cal T}}
\newcommand{\cJ}{{\cal J}}
\newcommand{\cS}{{\cal S}}
\newcommand{\cC}{{\cal C}}

\IEEEoverridecommandlockouts

\begin{document}
\title{An Efficient Primal-Dual Approach to Chance-Constrained Economic Dispatch}

\author{\IEEEauthorblockN{Gabriela Martinez\IEEEauthorrefmark{1},
Yu Zhang\IEEEauthorrefmark{2},
Georgios B. Giannakis\IEEEauthorrefmark{2}}
\IEEEauthorblockA{\IEEEauthorrefmark{1}Dept. of Biological \& Environmental Engineering,
Cornell University, Ithaca, USA \\
Email: gabriela.martinez@cornell.edu}
\IEEEauthorblockA{\IEEEauthorrefmark{2}Dept. of ECE and the Digital Technology Center,
University of Minnesota, Minneapolis, USA\\
Emails: \{zhan1220,georgios\}@umn.edu}
\thanks{* Supported by Consortium for Electric Reliability Technology Solutions (CERTS) grant PSERC/R\&M Project 3D, Cornell University.\newline \indent
$\dag$ Supported by Institute of Renewable Energy and the Environment (IREE) grant RL-0010-13, Univ. of Minnesota, and NSF-ECCS grant 1202135.}}

\maketitle

% -------------------------------------------------------------------------
% Abstract
% -------------------------------------------------------------------------
\begin{abstract}
To effectively enhance the integration of distributed and renewable energy sources in future smart microgrids,
economical energy management accounting for the principal challenge of the variable and non-dispatchable
renewables is indispensable and of significant importance.
Day-ahead economic generation dispatch with demand-side management for a microgrid in islanded mode is considered
in this paper. With the goal of limiting the risk of the loss-of-load probability,
a joint chance constrained optimization problem is formulated for the optimal multi-period energy scheduling with multiple wind farms.
Bypassing the intractable spatio-temporal joint distribution of the wind power
generation, a primal-dual approach is used  to obtain a suboptimal solution efficiently.
The method is based on first-order optimality conditions and successive approximation of
the probabilistic constraint by generation of $p$-efficient points.
Numerical results are reported to corroborate the merits of this approach.
\end{abstract}

\begin{IEEEkeywords}
Microgrids, renewable energy, economic dispatch, chance constraints, primal-dual approach.
\end{IEEEkeywords}

\section{Introduction}

\begin{figure}
\centering
\includegraphics[scale =0.62]{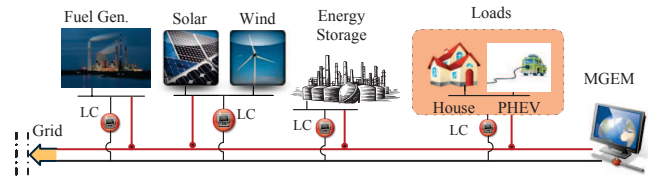}
\caption{The typical infrastructure of a microgrid
with communications (black) and energy flow (red) networks.
Each of the DERs and loads has a local controller (LC), which
coordinates with the microgrid energy manager
(MGEM) the scheduling of resources through the communications infrastructure.
}\label{fig:MGModel}
\end{figure}

As modern, small-scaled counterparts of the bulk power grids, microgrids are promising to achieve the goals of
improving efficiency, sustainability, security, and reliability of future electricity networks.
The very motivation of the infrastructure of microgrids is to bring power generation closer to the point where it is
consumed. In this way, distributed energy resources (DERs) and industrial, commercial, or residential electricity end-users
are deployed across a limited geographic area~\cite{Hatziargyriou-PESMag},
thereby incurring fewer thermal losses while bypassing other limitations imposed by the congested transmission networks.
Microgrids can operate either in grid-connected or islanded mode (a.k.a. isolated mode).
A typical microgrid configuration is depicted in Fig.~\ref{fig:MGModel}.
The communications between each local controller (LC) of DERs and controllable loads is coordinated via the microgrid energy manager (MGEM).

Besides the distributed storage (DS), a critical component comprising the DERs is
the renewable energy sources (RES) such as wind power, solar photovoltaic, biomass, hydro, and geothermal.
Because of the eco-friendly and price-competitive advantages,
renewable energy has been developing rapidly over the last few decades.
Recently, both the U.S. Department of Energy (DoE) and European Union (EU) proposed a very clear blueprint of
meeting $20$\% of the electricity consumption with renewables by 2030 and 2020, respectively.
Moreover, EU expects to achieve $80$-$95$\% greenhouse gas reductions by $100$\% renewable-integrated power systems by 2050~\cite{DOE08,EUreport}.

Accounting for the principal challenge incurred by the variability and uncertainty of the RES,
economical energy management plays an indispensable role toward achieving the goal of high penetration renewables.
Penalizing over- and under-estimation of wind energy, day-ahead economic dispatch (ED) is investigated in~\cite{HetzerYB08, YuNGGG-TSE13}.
A model predictive control (MPC)-based dynamic scheduling framework for variable wind generation and battery energy storage systems is proposed in~\cite{XieMPC12}.
Limiting the risk of the loss-of-load probability (LOLP), chance constrained optimization problems are formulated
for energy management problems including ED, unit commitment, and optimal power flow in~\cite{LiuX10, Qin13, FangGW12, BiChHa12, SjGaTo12}.
Globally optimal solutions are hard to obtain for the general non-convex chance-constrained problems.
Leveraging the Monte Carlo sampling based scenario approximation technique,
an efficient scalable solver is developed for chance constrained ED with correlated wind farms in~\cite{YuNGGG-ISGT13}.
However, it turns out to be too conservative in terms of the objective value, especially for the case of
a large-dimensional problem with a small prescribed risk requiring a large number of samples.

Building upon the earlier work in~\cite{YuNGGG-ISGT13},
the present paper proposes an efficient primal-dual approach to the multi-period ED and demand-side management (DSM) with multiple correlated wind farms.
The day-ahead economical scheduling task is formulated as a joint chance constrained optimization problem limiting the LOLP risk for an islanded microgrid.
The resulting non-convex chance constrained formulation is convexified using a convex combination of $p$-efficient points.
This procedure allows to decompose the generation of $p$-efficient points and approximate the convex hull of the original problem using these points.
The subproblem used to generate $p$-efficient points is equivalently rewritten as a mixed integer problem using a scenario approximation technique.
Comparing with the scenario approximation method proposed in~\cite{YuNGGG-ISGT13}, optimal solutions with smaller microgrid net costs are obtained by the novel approach
for different microgrids configurations. Numerical tests are implemented to corroborate the merits of the proposed approach.

The remainder of the paper is organized as follows.
Section~\ref{sec:form} formulates the chance constrained energy management problem, followed by
the development of the primal-dual approach in Section~\ref{sec:pridual}.
Numerical results are reported in Section~\ref{sec:sim}.
Finally, conclusions and research directions are provided in Section~\ref{sec:sum}.

\noindent \emph{Notation}. Boldface lower case letters represent vectors;
$\R^{n}$ and $\R_{+}^{n}$ are the space of $n \times 1$ vectors and the non-negative orthant, respectively;
$(\cdot)^{\top}$ stands for vector transpose;
$\mathbf{a} \succeq \mathbf{b}$ denotes the entry-wise inequality between two vectors; and
the probability of an event $A$ is denoted as $\Prob(A)$.

\section{Chance-constrained Economic Dispatch Formulation}\label{sec:form}

Consider an islanded microgrid featuring $M$ conventional  generators, $N$ dispatchable loads, $J$ storage units, and $I$ wind farms.
Let $\cT=\{1,2,\ldots,T\}$ be the scheduling horizon.
Denote $\cM = \{1,\dots,M\}$, $\cN=\{1,\dots,N\}$, and $\cJ = \{1,\dots,J\}$ as the index set of the conventional generation units,
the dispatchable loads, and storage units, respectively. Let $\pg{m}{t}$ be the power produced by the $m$th conventional generator,
$\pd{n}{t}$ the power consumed by the $n$th dispatchable load, $\pb{j}{t}$ the power delivered to the $j$th storage unit,
and $\bj{j}{t}$ the state of charge (SoC) of the $j$th storage unit at time $t$.
In addition to dispatchable loads, there is a fixed base load demand, $L^t$, that has to be served at each period.
The random wind power generated by the $i$th wind farm at time $t$ is denoted by $W_i^t,~i \in \cI := \{1,\dots,I\}$.
Let $W^t$ be the aggregated wind power at time $t$ defined as $W^t:=\sum_{i\in\cI}W_i^t$.

The goal of the day-ahead risk-averse ED is to minimize the microgrid operating cost
while satisfying the power demand with a prescribed high probability $p$.
Upon defining 
\begin{eqnarray*}
\bp_G &:=& [P_{G_1}^1,\dots,P_{G_1}^T,\dots,P_{G_M}^1,\dots,P_{G_M}^T] \\
\bp_D &:=& [P_{D_1}^1,\dots,P_{D_1}^T,\dots,P_{D_N}^1,\dots,P_{D_N}^T] \\
\bb &:=& [B_{1}^1,\dots,B_{1}^T,\dots,B_{J}^1,\dots,B_{J}^T] \\
\bp_B &:=& [P_{B_1}^1,\dots,P_{B_1}^T,\dots,P_{B_J}^1,\dots,P_{B_J}^T] \\
\end{eqnarray*}
the microgrid operating cost is given as:
\begin{align}
\label{eqn:obj}
&F(\bp_G,\bp_D,\bp_B,\bb) := \nonumber \\
&\sum_{t\in\cT}\Big(\sum_{m\in\cM}C_m(\pg{m}{t}) - \sum_{n\in\cN} U_n(\pd{n}{t}) + \sum_{j\in\cJ}H^t_j(\bj{j}{t})\Big)
\end{align}
where the function $C_m(\pg{m}{t})$ represents the $m$th conventional generation cost and it is assumed to be a convex quadratic function.
The function $U_n(\pd{n}{t})$ is a concave quadratic utility function of the $n$th dispatchable load.
The function $H^t_j(\bj{j}{t})$ is the $j$th storage usage cost, which is assumed to
be decreasing with respect to the state of charge of the $j$th storage unit (see also~\cite{YuNGGG-TSE13}).

Conventional generation constraints listed next represent (\ref{pgbounds}) power limit bounds,
(\ref{pgrup})-(\ref{pgrdn}) ramping up and down constraints, (\ref{pgsr}) spinning reserve.
Constraints (\ref{pdbound}) are the consumption bounds of the dispatchable loads.
Additional constraints are (\ref{sbound})-(\ref{sboundf}) SoC limits, (\ref{pbbound}) charging or discharging bounds,
(\ref{sdc}) SoC dynamic equations, and~\eqref{seff} storage efficiency constraints.
\begin{eqnarray}
&\pgmin{m}\leq \pg{m}{t}\leq\pgmax{m} & m\in\cM,t\in\cT \label{pgbounds}\\
&\pg{m}{t}-\pg{m}{t-1}\leq \rup{m} &m\in\cM,t\in\cT \label{pgrup}\\
&\pg{m}{t-1}-\pg{m}{t}\leq \rdn{m} &m\in\cM,t\in\cT \label{pgrdn}\\
&\sum_{m\in\cM}(\pgmax{m} - \pg{m}{t}) \geq SR^t &t\in\cT \label{pgsr}\\
&\pdmin{n}\leq \pd{}{}\leq\pdmax{n} &n\in\cN,t\in\cT \label{pdbound}\\
& \bj{j}{t}\leq\bjmax{j}&j\in\cJ,t\in\cT \label{sbound}\\
&\bj{j}{T}\geq\bjmin{j}&j\in\cJ \label{sboundf}\\
&{\pbmin{j}\leq\pb{j}{t}\leq\pbmax{j}}&j\in\cJ,t\in\cT\label{pbbound}\\
& \bj{j}{t}=\bj{j}{t-1}+\pb{j}{t}&j\in\cJ,t\in\cT \label{sdc}\\
& -\eta_j\bj{j}{t-1} \leq \pb{j}{t}&j\in\cJ,t\in\cT \label{seff}.
\end{eqnarray}

In addition to the aforementioned constraints, the most important
operation requirement is the power supply-demand balance.
However, for an islanded microgrid, energy transaction between the microgrid and the main grid is not applicable
for the supply-demand balance. In this case, a straightforward but effective way for
the day-ahead dispatch is to limit the risk of the LOLP, which is a metric
evaluating how frequent the total power supply can not satisfy the total demand.
Let $p\in (0,1)$ denote a prescribed probability level.
Restricting the joint LOLP can be equivalently written as:
\begin{equation}
\label{eqn:cc}
\Prob\big(L^t+\sum_{n\in\cN}\pd{n}{t}+\sum_{j\in\cJ}\pb{j}{t}-\sum_{m\in\cM}\pg{m}{t}\leq W^t,t\in\cT\big)\geq p.
\end{equation}

To this end, the risk-limiting ED task is tantamount to solving the following problem:
\begin{align*}
\label{ccprob}
\text{(P1)}\quad \min\, F(\bp_G,\bp_D,\bp_B,\bb),\,\, \text{subject to:}\,\, \eqref{pgbounds}-\eqref{eqn:cc}.
\end{align*}
Clearly, problem (P1) has a convex objective~\eqref{eqn:obj} as well as the linear equality and inequality constraints~\eqref{pgbounds}-\eqref{seff}.
Hence, the difficulty of solving (P1) lies in the joint chance constraint~\eqref{eqn:cc}.
The closed form of \eqref{eqn:cc} is intractable since the joint spatio-temporal
distribution of the wind power is unknown, and it is generally non-convex.

In the next section, in order to convexify the constraint~\eqref{eqn:cc} appropriately,
the definition of a $p$-efficient point will be introduced.
Interested readers are referred to \cite{shapiro2009lectures} for a comprehensive treatment of chance constrained problems and the corresponding optimality conditions.

\section{$p$-efficient points and Primal-Dual Method}\label{sec:pridual}

Let the random vector $\bw := [W^1,\dots,W^T]$ be the aggregated wind power outputs across the time slots. Let $\mathcal{Z}_p := \{\mathbf{z}\in \R^T : \Prob(\bw \succeq \mathbf{z})\geq p\}$. The $p$-efficient point is defined as follows~\cite[Sec.~4.3]{shapiro2009lectures}
\begin{definition}
Let $p\in (0,1)$. The vector $\be \in \R^T$ is called a $p$-efficient point if $\be\in\mathcal{Z}_p$ and
there is no $\bz\in\mathcal{Z}_p$ such that $\bz\succeq \be,\, \bz \ne \be$.
\end{definition}
Thus, a $p$-efficient point $\be\in\mathcal{Z}_p$ has maximal coordinates in $\mathcal{Z}_p$.  The convexification of constraint \eqref{eqn:cc} can be obtained as follows:
\begin{equation}
\label{eqn:ccp}
L^t+\sum_{n\in\cN}\pd{n}{t}+\sum_{j\in\cJ}\pb{j}{t}-\sum_{m\in \cM}\pg{m}{t}\leq u^t,\, t\in\cT
\end{equation}
where $\bu :=[u^1,\ldots,u^T]$ is a convex combination of $p$-efficient points.
The convexification of (P1) is thereby obtained by replacing~\eqref{eqn:cc} with~\eqref{eqn:ccp} as follows:
\begin{align*}
\label{ccprob-approx}
\text{(P2)}\quad \min\, F(\bp_G,\bp_D,\bp_B,\bb),\,\, \text{subject to:}\,\, \eqref{pgbounds}-\eqref{seff},\eqref{eqn:ccp}.
\end{align*}
%Let $\bm{\lambda} \in\R_{+}^T$ be the Lagrange multiplier associated with constraint (\ref{eqn:ccp}).
%Then, the dual functional associated with (P2) can be decomposed as:
%$\phi(\bm{\lambda}) = \phi_\rho(\bm{\lambda}) - \phi_\varsigma(\bm{\lambda})$.
%The function $\phi_\rho(\bm{\lambda})$ is defined as a deterministic minimization problem with respect to the primal variables $\{\bp_G,\bp_D,\bb,\bp_B\}$,
%while the stochastic dual subproblem $\phi_\varsigma(\bm{\lambda})$ is defined as follows:
%\begin{equation}
%\label{eqn:ppoint}
%\max_{\bv}\, \bm{\lambda}^\top \bv,\quad \text{ subject to:}\,\, \Prob(\bw \succeq \bv)\geq p.
%\end{equation}
%A detailed analytical description of the stochastic dual subproblem (\ref{eqn:ppoint}) is provided in~\cite[Lemma 4.67]{shapiro2009lectures}.
Let $\mathbf{x}$ collects variables $\{\bp_G,\bp_D,\bp_B,\bb\}$; 
$\mathcal{X} := \{\mathbf{x} : \mathbf{x}\in \eqref{pgbounds}-\eqref{seff}\}$; 
and $\bv := [v^1,\dots,v^T]$. Define further  
$\mathbf{g}(\mathbf{x}) := [g^1(\mathbf{x}),\dots,g^T(\mathbf{x})]$, where 
$g^t(\mathbf{x}):= L^t+\sum_{n\in\cN}\pd{n}{t}+\sum_{j\in\cJ}\pb{j}{t}-\sum_{m\in \cM}\pg{m}{t}, \forall t\in \cT$.
By splitting variables, (P1) can be equivalently reformulated as
\begin{align}
\min\, &F(\mathbf{x}) \\
\text{subject to:}\,&\mathbf{g}(\mathbf{x}) \preceq \mathbf{v} \label{ccsplt}\\
&{\mathbf{x}\in \mathcal{X}, \bv\in \mathcal{Z}_p}.
\end{align}
Let $\bm{\lambda} \in\R_{+}^T$ be the Lagrange multiplier associated with constraint \eqref{ccsplt}.
The partial Lagrangian function has the form
\begin{align}
\label{ccprob-Lagr}
L(\mathbf{x},\bv,\bm{\lambda}) =  F(\mathbf{x}) +  \bm{\lambda}^{\top}(\mathbf{g}(\mathbf{x})-\bv).
\end{align}
The dual function is thus obtained as
\begin{align}
\label{ccprob-dual}
\phi(\bm{\lambda}) = \inf_{\mathbf{x}\in \mathcal{X}, \bv \in \mathcal{Z}_p}\,
L(\mathbf{x},\bv,\bm{\lambda}) = \phi_\rho(\bm{\lambda}) - \phi_\varsigma(\bm{\lambda}),
\end{align}
where
\begin{align}
\phi_\rho(\bm{\lambda}) &:=\inf\{F(\mathbf{x})+\bm{\lambda}^{\top}\mathbf{g}(\mathbf{x}):\mathbf{x}\in \mathcal{X}\}, \label{dual-sub1}\\
\phi_\varsigma(\bm{\lambda}) &:=  \sup\{\bm{\lambda}^{\top} \bv:\mathbf{v}\in \mathcal{Z}_p\}. \label{dual-sub2}
\end{align}

To this end, Monte Carlo samples are needed to approximate the unknown joint distribution of the wind power output.
%$\cS:=\{1,\dots,N_s\}$,
Let $\bw_{s}:=[W^1_s,\dots,W^T_s], s\in\cS:=\{1,\dots,N_s\}$ denote an independent and identically distributed (i.i.d.) sample of the random vector $\bw$.
Let $\1_{A}$ be the indicator function of an event $A$; i.e., $\1_{A}$ takes the value $1$ if $A$ is true, and the value $0$ otherwise.
The empirical survival function is defined as
%\begin{equation*}
$\cP(\bw)=\frac{1}{N_s}\sum_{s\in\cS}\1_{\{\bw_s\succeq \bw\}}.$
%\end{equation*}
Similarly, the marginal empirical survival function is given as $\cP_t(w)=\frac{1}{N_s}\sum_{s\in\cS}\1_{\{W^t_s\geq w\}}.$

The optimization problem~\eqref{dual-sub2} can be approximated with the following mixed integer problem:
\begin{subequations}
\label{eqn:mip}
\begin{align}
\max\, &\bm{\lambda}^\top \bv \\
\text{subject to:}\,\,
& \bv-\bl_p \preceq (\bw_s-\bl_p)z_s,\, s\in \cS\\
& \sum_{s\in\cS}z_s\geq pN_s\\
& z_s\in\{0,1\},s\in\cS.
\end{align}
\end{subequations}
where $\bl_p :=[\ell^1_p,\ldots,\ell^T_p] \in\R^T$ is such that for $t\in\cT$, $\cP_t(\ell_p^t)\geq p$. Approximations of the $p$-efficient points can be obtained by solving  problem~\eqref{eqn:mip}.
%(see details in~\cite{dbr,ddgm1,ddgm2,dpr}).
As tabulated in Algorithm~\ref{algo:prim-dual}, the proposed primal-dual method approximates
the feasible set of the optimization problem by generating $p$-efficient points at each iteration.
The collection of $p$-efficient points is then used to approximate constraint~\eqref{eqn:ccp}.

The solution found by the primal-dual method is $\epsilon$-optimal with respect to the sample for (P2) (see details in~\cite{dbr,ddgm1,ddgm2,dpr}). Let $\cK_a$ denote the set of active optimal $p$-efficient points: $\cK_{a}=\{k\in\cK:\alpha_k>0\}$. If $\cK_{a}$ contains only one element, then the solution found is $\epsilon$-optimal with respect to the sample for (P1). Otherwise, the optimal value found by the method is a lower bound for the optimal value of (P1). Therefore, primal feasible points for (P1) can be obtained by solving for $k\in\cK_{a}$ the following optimization problem:
\begin{subequations}
\begin{align}
&\min\,\, F(\bp_G,\bp_D,\bp_B,\bb)\\
&\text{subject to:}\quad (\ref{pgbounds})-(\ref{seff})\\
&L^t+\sum_{n\in\cN}\pd{n}{t}+\sum_{j\in\cJ}\pb{j}{t}-\sum_{m\in\cM}\pg{m}{t}\leq v_k^t.
\end{align}
\end{subequations}

%%%%%%%%%%%%%%%%%%%%%%%%%%%%%%%%%%%%%%%%%%%%%%%%%%%%
\begin{algorithm}[t]
\caption{Primal-dual approach to chance-constrained ED}
\label{algo:prim-dual}
\begin{algorithmic}[1]
\State Let $\epsilon >0$, $k=1$, $\cK=\{k\}$, $\bm{\lambda}^k\in\R_{+}^T$. Find $\bv_k$ solution of (\ref{eqn:mip}) with $\bm{\lambda}^k$.
\State Solve the master problem:
\begin{align}
&\min F(\bp_G,\bp_D,\bp_B,\bb) \notag \\
&\text{subject to:}\quad (\ref{pgbounds})-(\ref{seff}) \notag \\
&L^t+\sum_{n\in\cN}\pd{n}{t}+\sum_{j\in\cJ}\pb{j}{t}-\sum_{m\in\cM}\pg{m}{t}\leq \sum_{k\in\cK}\alpha_kv_k^t \label{eqn:approx}\\
&\sum_{k\in\cK}\alpha_k = 1, \alpha_k\geq 0, k\in\cK. \notag
\end{align}
Let $\bm{\lambda}^{k+1}$ be the multiplier associated with constraint~\eqref{eqn:approx}.
\State Solve the dual subproblem (\ref{eqn:mip}) with $\bm{\lambda}^{k+1}$. Let $\bar{\phi}_\varsigma$ be its optimal value and $\bv_{k+1}$ its solution.
\State Define $\phi_k := \max_{j\in\cK}(\bm{\lambda}^{k+1})^\top \bv_j$. If $|\bar{\phi}_\varsigma-\phi_k|<\epsilon$ stop. Otherwise $k\rightarrow k+1$, $\cK=\cK\cup\{k\}$, and go to step 2.
\end{algorithmic}
\end{algorithm}
%%%%%%%%%%%%%%%%%%%%%%%%%%%%%%%%%%%%%%%%%%%%%%%%%%%%%%

Two remarks are in order on the $p$-efficient points and the optimal values.
\begin{remark}
\label{rmk:pmin}
Consider the sample $\{\bw_s, s\in\cS\}$. Let $p\in (0,1)$, $\be$ be a $p$-efficient point in the sample, and $\bar{\bw}:=[\bar{w}^1,\ldots,\bar{w}^T]$ be the minimal value in the sample.
Define $\cC_p :=\{\bw\in\R^T, \cP(\bw)\geq p\}$, the $p$-efficient points are thus on the boundary of the set $\cC_p$; see \cite[Theorem 4.60]{shapiro2009lectures}.
Furthermore, $\bar{\bw}$ is in the interior of the set $\cC_p$ since $\cP(\bar{\bw}) = 1\geq p$.
It follows from the definition of $p$-efficient point that $\bar{\bw}\ne\be$, $\bar{\bw} \preceq \be$.
In addition, let $\bu$ be a convex combination of $p$-efficient points then $\bar{\bw}\preceq\bu$.
\end{remark}
\begin{remark}
It follows from Remark~\ref{rmk:pmin} that any feasible solution to the problem:
\begin{subequations}
\label{eqn:robust}
\begin{align}
&\min\,\, F(\bp_G,\bp_D,\bp_B,\bb)\\
&\text{subject to:}\quad (\ref{pgbounds})-(\ref{seff})\\
&L^t+\sum_{n\in\cN}\pd{n}{t}+\sum_{j\in\cJ}\pb{j}{t}-\sum_{m\in\cM}\pg{m}{t}\leq \bar{w}^t.
\end{align}
\end{subequations}
is feasible for the optimization problems (P1) and (P2).
%\end{remark}
%\begin{remark}
%\label{rmk:values}
Denote by $F_{\be}$ and $F_{\bu}$ the optimal values of the optimization problems  (P1) and (P2), respectively.
Let $F_{\bar{\bw}}$ be the optimal value of (\ref{eqn:robust}). Then, $F_{\bu}\leq F_{\be}\leq F_{\bar{\bw}}$.
\end{remark}

\begin{table}[t]
\centering
\caption{Generation limits, ramping rates, and cost coefficients.
The units of $a_m$ and $b_m$ are \$/(kWh)$^{2}$ and \$/kWh, respectively.}\label{tab:generator}
    \begin{tabular}{  c | c  c | c  c | c  c }
    \hline
Unit &$P_{G_{\text{min},m}}$ &$P_{G_{\text{max},m}}$  &$R_{\textrm{up},m}$ &$R_{\textrm{dn},m}$  &$a_m$           &$b_m$       \\ \hline
1    & 10                     & 30            & 15                & 15               &0.006           & 0.5                      \\
2    & 8                      & 50             & 40                 & 40               &0.003           & 0.25                    \\
3    & 15                     & 70             & 20                 & 20               &0.004           & 0.3                         \\
    \hline
    \end{tabular}
\end{table}

\begin{table}[t]
\centering
\caption{Parameters of dispatchable loads.
The units of $c_n$ and $d_n$ are \$/(kWh)$^{2}$ and \$/kWh, respectively.}\label{tab:load}
    \begin{tabular}{ c | c  c  c  c  c  c }
    \hline
    Load          &  1          &   2          &   3          &  4       &  5         &  6   \\  \hline
    $P_{D_{\text{min},n}}$           & 1.5        & 3.3      & 2        & 5.7       & 4        & 9     \\
    $P_{D_{\text{max},n}}$           & 8          & 10       & 15        & 24       & 20        & 35   \\ \hline
    $c_n\times 10^3$                    &-4.5   &-1.1   &-1.9   &-1.3   &-1.4   &-2.6    \\
    $d_n$                       &0.15    &0.37    &0.62    &0.44    &0.45   &0.87  \\
\hline
    \end{tabular}
\end{table}

% ----------------------------
\section{Numerical Tests}
\label{sec:sim}
% ----------------------------

Numerical case studies are provided to verify the effectiveness of the novel approach in this section.
The solver~\texttt{Cplex 12.5} and the Python-based modeling package~\texttt{Pyomo} are used to solve the subproblem listed in the primal-dual algorithm.
The tested microgrid in island mode consists of $M=3$ conventional generators, $N=6$ dispatchable loads, $J=3$ storage units, and $I=4$ wind farms.
The scheduling horizon spans $T=8$ hours, corresponding to the interval $4$pm--$12$am.
The generation costs $C_m(P_{G_m})=a_mP_{G_m}^2+b_mP_{G_m}$ and the load utilities $U_n(P_{D_n})=c_nP_{D_n}^2+d_nP_{D_n}$
are set to be time-invariant, of which the  parameters are listed in Tables~\ref{tab:generator} and~\ref{tab:load}.
The spinning reserve is set to be zero;
The storage  usage cost $H^t_j(\bj{j}{t}):= \beta_j^t(\bjmax{j}-\bj{j}{t})$ is set to be linear.
The state of charge bounds are $\bjmin{j}=5$, $\bjmax{j}=30$, for $j\in \cJ$.
The storage usage weights used are $\beta_1^t=0.05t$, $\beta_2^t=\beta_3^t=0.1$.
The fixed base load demand used is  $L^t =[43.35, 43.95, 48, 48.825, 46.125, 44.1, 41.625, 38.25]$ (kWh).

To estimate the required $p$-efficient points, Monte Carlo samples are obtained by a turbine-specific
wind energy conversion system (WECS), where the wind speed samples are generated using a two-parameter $(c,k)$ Weibull distribution.
The relevant parameters of the WECS are listed in Tables~\ref{tab:wecs}.
An autoregressive model is utilized to capture the temporal correlation of the wind speed across the horizon.
The lag-one temporal correlations are chosen as $\{\phi_i\}_{i=1}^{I}:=\{0.15,0.43,0.67,0.59\}$; and
the spatial correlation coefficient matrix is given as (see also~\cite{YuNGGG-ISGT13})
\begin{align*}
\mathbf{C}=
\left[\begin{array}{cccc}
1    &0.1432    &0.4388   &-0.0455\\
0.1432    &1   &-0.4555    &0.8097\\
0.4388   &-0.4555    &1   &-0.7492\\
-0.0455    &0.8097   &-0.7492    &1
\end{array}\right].
\end{align*}

\begin{table}[t]
\centering
\caption{Parameters of the WECS}\label{tab:wecs}
    \begin{tabular}{ c | c | c | c | c } \hline
Parameter     &$c$          &$k$        &$v_{\textrm{in,rated,out}}$ (m/s)       &$w_{\textrm{rated}}$ (kWh)  \\  \hline
Value          &10           &2.2              &3,14,26                               &10 \\
\hline
\end{tabular}
\end{table}

\begin{table}[t]
\centering
\caption{Optimal microgrid net costs of the novel method and the SAA in~\cite{YuNGGG-ISGT13}.}
\label{tab:cost}
    \begin{tabular}{ c | c | c | c || c }
    \hline
                    & $p=0.9$          & $p=0.95$            & $p=0.99$   & SAA~\cite{YuNGGG-ISGT13}   \\  \hline
    $N_s = 100$        &54.67  &60.59    &67.09     & 68.86  \\
    $N_s = 500$         &63.64   &69.94   &76.98     & 81.21           \\
    $N_s = 1000$         &63.16   &68.82   &77.77    & 82.84            \\
 \hline
    \end{tabular}
\end{table}

\begin{figure}[t]
\centering
\includegraphics[scale=.5]{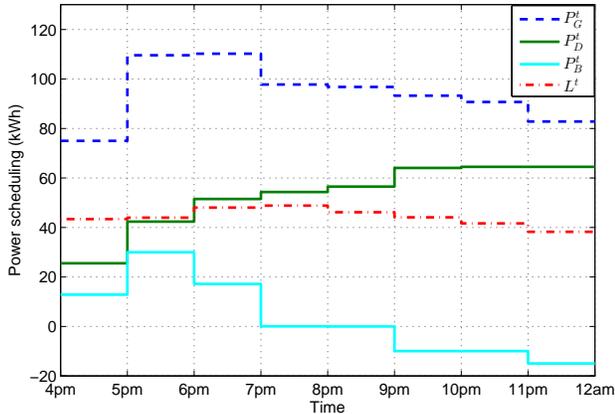}
%\vspace{-0.2cm}
\caption{Optimal power schedules.}
\label{fig:sche}
%\vspace{-0.4cm}
\end{figure}

\begin{figure}[t]
\centering
\includegraphics[scale=.5]{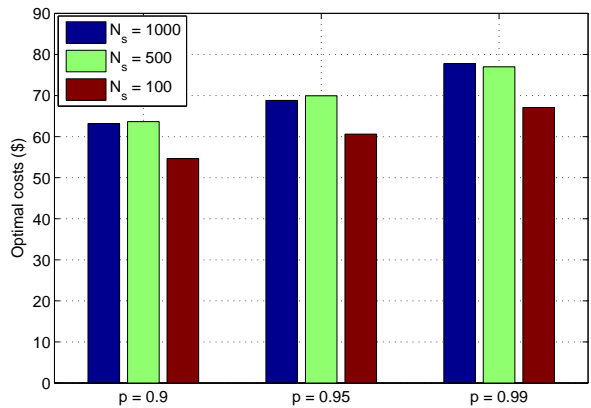}
\caption{Optimal microgrid net costs.}
\label{fig:cost}
\end{figure}

%\begin{table}[th]
%\centering
%\caption{Prescribed LOLP $\alpha_{pr}$ versus actual LOLP $\alpha_{ac}$.}\label{tab:LOLPvalidation}
%    \begin{tabular}{ c | c | c | c  }
%    \hline
%    $\alpha_{pr}$          &0.9    &0.95    &0.99    \\  \hline
%    $\alpha_{ac}$          &0.0002    &0.0346    &0.0464     \\
%\hline
%    \end{tabular}
%\end{table}

The optimal power schedules are depicted in Fig.~\ref{fig:sche} for the case of $p=0.95$ and $N_s=1,000$.
The stair step curves include $P_G^t:=\sum_mP_{G_m}^t$, $P_D^t :=\sum_nP_{D_n}^t$, and $P_B^t :=\sum_j P_{B_j}^t$
denoting the total conventional power, total elastic demand, and total (dis)charging power, respectively.
As expected, the conventional power $P_G^t$ exhibits a similar trend with the fixed base load demand $L^t$
across the time horizon. In addition, the elastic demand $P_D^t$ trends oppositely to $L^t$ reflecting
the peak load shifting ability of the proposed scheme. Specifically,
by smartly scheduling the high power demand of the deferred load $P_D^t$ to the slots $9$pm-$12$am,
the potential peak of the total demand from $6$pm to $9$pm can be shaved.
Recall that storage units play a role of loads (generators) whenever $P_B^t \geq 0$ ($P_B^t < 0$)
because they are charging (discharging) power from (to) the microgrid. As shown in Fig.~\ref{fig:sche},
it can be seen that from $4$pm-$7$pm, storage units are charging since the total power demand of $P_D^t$ and $L^t$ are relatively low.
From $7$pm-$12$am, storage units keep discharging to generate energy supporting the dispatchable and base loads.
Consequently, this reduces the conventional generation during these time horizons.

Figure~\ref{fig:cost} shows the optimal microgrid net costs for different values of the probability level $p$
and numbers of i.i.d. samples. Clearly, the optimal cost becomes higher with the increase of the $p$ values for
each of three different number of samples.
This is because a less operation risk is allowed for a larger $p$ value, resulting in a more restricting feasible set
of the primal variables, i.e., the power schedules $\{\bp_G,\bp_D,\bb \}$, and hence a higher net cost.
Interestingly, it can be seen that for a given $p$ probability, the net costs corresponding to
$N_s = 500$ and $N_s = 1000$ are roughly the same. In fact, the estimation of $p$-efficient points
becomes more accurate with more number of samples. Therefore, the solution obtained by the
$p$-efficient points based primal-dual method is not sensitive to the number of samples
if it is already large enough to estimate the efficient points accurately.
It is worth pointing out that this is not the case for the scenario approximation approach (SAA) proposed in~\cite{YuNGGG-ISGT13}.
The effective wind power resource given by the SAA essentially boils down to the
worst case scenario $\bar{\bw}:=\min_{s\in \cS}\{{\bw}^s\}$ which may be decreasing with the increasing size of samples.

Finally, the optimal microgrid costs of the novel method and the SAA~\cite{YuNGGG-ISGT13} are listed in Table~\ref{tab:cost}.
Note that for a fixed number of samples $N_s$, the optimal costs
obtained by the SAA are identical for different $p$ values.
As detailed in Remark 2, the optimal microgrid net costs of the proposed method are consistently lower than those
obtained by the SAA for different values of $N_s$ and $p$. This fact shows the
advantage of the novel primal-dual approach overcoming the potential conservativeness of the SAA, and
capability of obtaining more economical operation points for the day-ahead microgrid power dispatch.

\section{Conclusions and Future Work}\label{sec:sum}

In this paper, the day-ahead ED with DSM for islanded microgrids with spatio-temporal wind farms is considered.
The power scheduling task is formulated as a chance constrained optimization problem based on the LOLP.
Leveraging $p$-efficient points, a primal-dual approach is developed for efficiently solving the resulting non-convex
chance constrained problem. Case studies corroborate the effectiveness of the proposed approach that is capable of
finding economical power operation points.

Some interesting research directions are worthy of exploring towards improving the work presented in this paper.
When the microgrid is operated in islanded mode, changes in power demand can cause changes in frequency and voltage levels. Therefore, frequency regulation is important to maintain system stability. The economic dispatch model proposed in this paper could be improved by incorporating a frequency control method, for example, in \cite{freqregwind} wind turbine controls are used to regulate the frequency of the microgrid, other frequency control methods are reviewed in \cite{freqreg}. A realistic modeling of ancillary services, as the multi-agent model proposed in \cite{multiagent}, could be used to determine adequate spinning reserve levels for system. In \cite{storagesandia} are proposed accurate models for energy storage devices. Since energy storage is an important component for an islanded microgrid, the constraints modeling the storage devices could be improved by considering the properties detailed in \cite{storagesandia}.  

The mixed integer formulation (\ref{eqn:mip}) has a knapsack constraint which is known to be NP-hard.
Strength formulations as in~\cite{Luedtke2007,Bienstock01} may be useful to generate the $p$-efficient points with reduced computational complexity. The solution methods proposed in~\cite{ddgm1} with a strength formulation of~\eqref{eqn:mip}
can be used to manage large size samples.

%\IEEEtriggercmd{\enlargethispage{-0.52in}}
%\IEEEtriggeratref{1}

%%%%%%%%%%%%%%%%%%%%%%%%%%%%%%%%%%%%%%%%%%%%%%
\bibliographystyle{IEEEtran}
\bibliography{biblio}

\begin{thebibliography}{10}
\providecommand{\url}[1]{#1}
\csname url@samestyle\endcsname
\providecommand{\newblock}{\relax}
\providecommand{\bibinfo}[2]{#2}
\providecommand{\BIBentrySTDinterwordspacing}{\spaceskip=0pt\relax}
\providecommand{\BIBentryALTinterwordstretchfactor}{4}
\providecommand{\BIBentryALTinterwordspacing}{\spaceskip=\fontdimen2\font plus
\BIBentryALTinterwordstretchfactor\fontdimen3\font minus
  \fontdimen4\font\relax}
\providecommand{\BIBforeignlanguage}[2]{{%
\expandafter\ifx\csname l@#1\endcsname\relax
\typeout{** WARNING: IEEEtran.bst: No hyphenation pattern has been}%
\typeout{** loaded for the language `#1'. Using the pattern for}%
\typeout{** the default language instead.}%
\else
\language=\csname l@#1\endcsname
\fi
#2}}
\providecommand{\BIBdecl}{\relax}
\BIBdecl

\bibitem{Hatziargyriou-PESMag}
N.~Hatziargyriou, H.~Asano, R.~Iravani, and C.~Marnay, ``Microgrids: An
  overview of ongoing research, development, and demonstration projects,''
  \emph{IEEE Power \& Energy Mag.}, vol.~5, no.~4, pp. 78--94, July--Aug. 2007.

\bibitem{DOE08}
``20\% wind energy by 2030: Increasing wind energy's contribution to {U.S.}
  electricity supply,'' July 2008, [Online]. Available:
  \url{http://www1.eere.energy.gov/wind/pdfs/41869.pdf}.

\bibitem{EUreport}
{European Wind Energy Association}, ``{EU} energy policy to 2050 -- achieving
  80-95\% emissions reductions,'' Tech. Rep., Mar. 2011.

\bibitem{HetzerYB08}
J.~Hetzer, C.~Yu, and K.~Bhattarai, ``An economic dispatch model incorporating
  wind power,'' \emph{IEEE Trans. on Energy Conver.}, vol.~23, no.~2, pp.
  603--611, Jun. 2008.

\bibitem{YuNGGG-TSE13}
Y.~Zhang, N.~Gatsis, and G.~B. Giannakis, ``Robust energy management for
  microgrids with high-penetration renewables,'' \emph{IEEE Trans. on
  Sustainable Energy}, vol.~4, no.~4, pp. 944--953, Oct. 2013.

\bibitem{XieMPC12}
L.~Xie, Y.~Gu, A.~Eskandari, and M.~Ehsani, ``Fast {MPC}-based coordination of
  wind power and battery energy storage systems,'' \emph{Journal of Energy
  Eng.}, vol. 138, no.~2, pp. 43--53, June 2012.

\bibitem{LiuX10}
X.~Liu and W.~Xu, ``Economic load dispatch constrained by wind power
  availability: A here-and-now approach,'' \emph{IEEE Trans. on Sustainable
  Energy}, vol.~1, no.~1, pp. 2--9, Apr. 2010.

\bibitem{Qin13}
J.~Qin, B.~Zhang, and R.~Rajagopal, ``Risk limiting dispatch with ramping
  constraints,'' in \emph{Proc. of IEEE Intl. Conf. on Smart Grid Commun.},
  Vancouver, Canada, Oct. 2013.

\bibitem{FangGW12}
Q.~Fang, Y.~Guan, and J.~Wang, ``A chance-constrained two-stage stochastic
  program for unit commitment with uncertain wind power output,'' \emph{IEEE
  Trans. on Power Syst.}, vol.~27, no.~1, pp. 206--215, 2012.

\bibitem{BiChHa12}
D.~Bienstock, M.~Chertkov, and S.~Harnett, ``Chance constrained optimal power
  flow: Risk-aware network control under uncertainty,'' Sept. 2012, [Online].
  Avaialble: \url{http://arxiv.org/pdf/1209.5779.pdf}.

\bibitem{SjGaTo12}
E.~Sj\"{o}din, D.~F. Gayme, and U.~Topcu, ``Risk-mitigated optimal power flow
  for wind powered grids,'' in \emph{Proc. American Control Conf.}, Montreal,
  Canada, June 2012, pp. 4431--4437.

\bibitem{YuNGGG-ISGT13}
Y.~Zhang, N.~Gatsis, and G.~B. Giannakis, ``Risk-constrained energy management
  with multiple wind farms,'' in \emph{Proc. of Innovative Smart Grid Tech.},
  Washington, D.C., Feb. 2013.

\bibitem{shapiro2009lectures}
A.~Shapiro, D.~Dentcheva, and A.~P. Ruszczy{\'n}ski, \emph{Lectures on
  stochastic programming: modeling and theory}.\hskip 1em plus 0.5em minus
  0.4em\relax SIAM, 2009, vol.~9.

\bibitem{dbr}
\BIBentryALTinterwordspacing
D.~Dentcheva, B.~Lai, and A.~Ruszczy\'{n}ski,
  ``\BIBforeignlanguage{English}{Dual methods for probabilistic optimization
  problems*},'' \emph{\BIBforeignlanguage{English}{Mathematical Methods of
  Operations Research}}, vol.~60, no.~2, pp. 331--346, 2004. [Online].
  Available: \url{http://dx.doi.org/10.1007/s001860400371}
\BIBentrySTDinterwordspacing

\bibitem{ddgm1}
\BIBentryALTinterwordspacing
D.~Dentcheva and G.~Martinez, ``\BIBforeignlanguage{English}{Regularization
  methods for optimization problems with probabilistic constraints},''
  \emph{\BIBforeignlanguage{English}{Mathematical Programming}}, vol. 138, no.
  1-2, pp. 223--251, 2013. [Online]. Available:
  \url{http://dx.doi.org/10.1007/s10107-012-0539-6}
\BIBentrySTDinterwordspacing

\bibitem{ddgm2}
\BIBentryALTinterwordspacing
------, ``\BIBforeignlanguage{English}{Augmented lagrangian method for
  probabilistic optimization},'' \emph{\BIBforeignlanguage{English}{Annals of
  Operations Research}}, vol. 200, no.~1, pp. 109--130, 2012. [Online].
  Available: \url{http://dx.doi.org/10.1007/s10479-011-0884-5}
\BIBentrySTDinterwordspacing

\bibitem{dpr}
\BIBentryALTinterwordspacing
D.~Dentcheva, A.~Prekopa, and A.~Ruszczynski,
  ``\BIBforeignlanguage{English}{Concavity and efficient points of discrete
  distributions in probabilistic programming},''
  \emph{\BIBforeignlanguage{English}{Mathematical Programming}}, vol.~89,
  no.~1, pp. 55--77, 2000. [Online]. Available:
  \url{http://dx.doi.org/10.1007/PL00011393}
\BIBentrySTDinterwordspacing

\bibitem{freqregwind}
K.~Bunker and W.~Weaver, ``Microgrid frequency regulation using wind turbine
  controls,'' in \emph{Power and Energy Conference at Illinois (PECI), 2014},
  Feb 2014, pp. 1--6.

\bibitem{freqreg}
S.~Baudoin, I.~Vechiu, and H.~Camblong, ``A review of voltage and frequency
  control strategies for islanded microgrid,'' in \emph{System Theory, Control
  and Computing (ICSTCC), 2012 16th International Conference on}, Oct 2012, pp.
  1--5.

\bibitem{multiagent}
N.-P. Yu, C.-C. Liu, and J.~Price, ``Evaluation of market rules using a
  multi-agent system method,'' \emph{Power Systems, IEEE Transactions on},
  vol.~25, no.~1, pp. 470--479, Feb 2010.

\bibitem{storagesandia}
J.~Eyer and G.~Corey, ``Energy storage for the electricity grid: Benefits and
  market potential assessment guide,'' {S}andia {N}ational {L}aboratories,
  Tech. Rep., [Online]. Available:
  \url{https://www.smartgrid.gov/sites/default/files/resources/energy_storage.%
pdf}.

\bibitem{Luedtke2007}
J.~Luedtke, S.~Ahmed, and G.~Nemhauser, ``An integer programming approach for
  linear programs with probabilistic constraints,'' in \emph{Proceedings of the
  12th International Conference on Integer Programming and Combinatorial
  Optimization}, ser. IPCO '07.\hskip 1em plus 0.5em minus 0.4em\relax
  Springer-Verlag, 2007, pp. 410--423.

\bibitem{Bienstock01}
D.~Bienstock, ``Approximate formulations for 0-1 knapsack sets,''
  \emph{Operations Research Letters}, vol.~36, no.~3, pp. 317 -- 320, 2008.

\end{thebibliography}
%%%%%%%%%%%%%%%%%%%%%%%%%%%%%%%%%%%%%%%%%%%%%%

\end{document}